# Asymptotic oracle properties of SCAD-penalized least squares estimators

Jian Huang[1] and Huiliang Xie[1]

*University of Iowa*

**Abstract:** We study the asymptotic properties of the SCAD-penalized least squares estimator in sparse, high-dimensional, linear regression models when the number of covariates may increase with the sample size. We are particularly interested in the use of this estimator for simultaneous variable selection and estimation. We show that under appropriate conditions, the SCAD-penalized least squares estimator is consistent for variable selection and that the estimators of nonzero coefficients have the same asymptotic distribution as they would have if the zero coefficients were known in advance. Simulation studies indicate that this estimator performs well in terms of variable selection and estimation.

## 1. Introduction

Consider a linear regression model

$$Y = \beta_0 + \mathbf{X}'\boldsymbol{\beta} + \varepsilon,$$

where $\boldsymbol{\beta}$ is a $p \times 1$ vector of regression coefficients associated with $\mathbf{X}$. We are interested in estimating $\boldsymbol{\beta}$ when $p \to \infty$ as the sample size $n \to \infty$ and when $\boldsymbol{\beta}$ is sparse, in the sense that many of its elements are zero. This is motivated from biomedical studies investigating the relationship between a phenotype of interest and genomic covariates such as microarray data. In many cases, it is reasonable to assume a sparse model, because the number of important covariates is usually relatively small, although the total number of covariates can be large.

We use the SCAD method to achieve variable selection and estimation of $\boldsymbol{\beta}$ simultaneously. The SCAD method is proposed by Fan and Li [1] in a general parametric framework for variable selection and efficient estimation. This method uses a specially designed penalty function, the smoothly clipped absolute deviation (hence the name SCAD). Compared to the classical variable selection methods such as subset selection, the SCAD has two advantages. First, the variable selection with SCAD is continuous and hence more stable than the subset selection, which is a discrete and non-continuous process. Second, the SCAD is computationally feasible for high-dimensional data. In contrast, computation in subset selection is combinatorial and not feasible when $p$ is large. In addition to the SCAD method, several other penalized methods have also been proposed to achieve variable selection and estimation simultaneously. Examples include the bridge penalty (Frank and Friedman [3]), LASSO (Tibshirani [11]), and the Elastic-Net (Enet) penalty (Zou and Hastie [14]), among others.

Fan and Li [1] and Fan and Peng [2] studied asymptotic properties of SCAD penalized likelihood methods. Their results are concerned with local maximizers of the penalized likelihood, but not the maximum penalized estimators. These results do not imply existence of an estimator with the properties of the local maximizer without auxiliary information about the true parameter value. Therefore, they are not applicable to the SCAD-penalized maximum likelihood estimators, nor the SCAD-penalized estimator. Knight and Fu [7] studied the asymptotic distributions of bridge estimators when the number of covariates is fixed. Huang, Horowitz and Ma [4] studied the bridge estimators with a divergent number of covariates in a linear regression model. They showed that the bridge estimators have an oracle property under appropriate conditions if the bridge index is strictly between 0 and 1. Several earlier studies have investigated the properties of regression estimators with a divergent number of covariates. See, for example, Huber [5] and Portnoy [9, 10]. Portnoy proved consistency and asymptotic normality of a class of M-estimators of regression parameters under appropriate conditions. However, he did not consider penalized regression or selection of variables in sparse models.

In this paper, we study the asymptotic properties of the SCAD-penalized least squares estimator, abbreviated as LS-SCAD estimator henceforth. We show that the LS-SCAD estimator can correctly select the nonzero coefficients with probability converging to one and that the estimators of the nonzero coefficients are asymptotically normal with the same means and covariances as they would have if the zero coefficients were known in advance. Thus, the LS-SCAD estimators have an oracle property in the sense of Fan and Li [1] and Fan and Peng [2]. In other words, this estimator is asymptotically as efficient as the ideal estimator assisted by an oracle who knows which coefficients are nonzero and which are zero.

The rest of this article is organized as follows. In Section 2, we define the LS-SCAD estimator. The main results for the LS-SCAD estimator are given in Section 3, including the consistency and oracle properties. Section 4 describes an algorithm for computing the LS-SCAD estimator and the criterion to choose the penalty parameter. Section 5 offers simulation studies that illustrate the finite sample behavior of this estimator. Some concluding remarks are given in Section 6. The proofs are relegated to the Appendix.

## 2. Penalized regression with the SCAD penalty

Let $(\mathbf{X}_i, Y_i), i = 1, \ldots, n$ be $n$ observations satisfying

$$Y_i = \beta_0 + \mathbf{X}'_i \boldsymbol{\beta} + \varepsilon_i, \quad i = 1, \ldots, n,$$

where $Y_i \in \mathrm{R}$ is a response variable, $\mathbf{X}_i$ is a $p_n \times 1$ covariate vector and $\varepsilon_i$ has mean 0 and variance $\sigma^2$. Here the superscripts are used to make it explicit that both the covariates and parameters may change with $n$. For simplicity, we assume $\beta_0 = 0$. Otherwise we can center the covariates and responses first.

In sparse models, the $p_n$ covariates can be classified into two categories: the important ones whose corresponding coefficients are nonzero and the trivial ones whose coefficients are zero. For convenience of notation, we write

$$\boldsymbol{\beta} = (\boldsymbol{\beta}'_1, \boldsymbol{\beta}'_2)',$$

where $\boldsymbol{\beta}'_1 = (\beta_1, \ldots, \beta_{k_n})$ and $\boldsymbol{\beta}'_2 = (0, \ldots, 0)$. Here $k_n (\leq p_n)$ is the number of nontrivial covariates. Let $m_n = p_n - k_n$ be the number of zero coefficients. Let

$\mathbf{Y} = (Y_1, \ldots, Y_n)'$ and let $\mathbb{X} = (X_{ij}, 1 \leq i \leq n, 1 \leq j \leq p_n)$ be the $n \times p_n$ design matrix. According to the partition of $\boldsymbol{\beta}$, write $\mathbb{X} = (\mathbb{X}_1, \mathbb{X}_2)$, where $\mathbb{X}_1$ and $\mathbb{X}_2$ are $n \times k_n$ and $n \times m_n$ matrices, respectively.

Given $a > 2$ and $\lambda > 0$, the SCAD penalty at $\theta$ is

$$p_\lambda(\theta; a) = \begin{cases} \lambda|\theta|, & |\theta| \leq \lambda, \\ -(\theta^2 - 2a\lambda|\theta| + \lambda^2)/[2(a-1)], & \lambda < |\theta| \leq a\lambda, \\ (a+1)\lambda^2/2, & |\theta| > a\lambda. \end{cases}$$

More insight into it can be gained through its first derivative:

$$p'_\lambda(\theta; a) = \begin{cases} \operatorname{sgn}(\theta)\lambda, & |\theta| \leq \lambda, \\ \operatorname{sgn}(\theta)(a\lambda - |\theta|)/(a-1), & \lambda < |\theta| \leq a\lambda, \\ 0, & |\theta| > a\lambda. \end{cases}$$

The SCAD penalty is continuously differentiable on $(-\infty, 0) \cup (0, \infty)$, but not differentiable at 0. Its derivative vanishes outside $[-a\lambda, a\lambda]$. As a consequence, SCAD penalized regression can produce sparse solutions and unbiased estimates for large coefficients. More detailed discussions of this penalty can be found in Fan and Li (2001).

The penalized least squares objective function for estimating $\boldsymbol{\beta}$ with the SCAD penalty is

$$(1) \qquad Q_n(\mathbf{b}; \lambda_n, a) = \|\mathbf{Y} - \mathbb{X}\mathbf{b}\|^2 + n \sum_{j=1}^{p_n} p_{\lambda_n}(b_j; a),$$

where $\|\cdot\|$ is the $L_2$ norm. Given penalty parameters $\lambda_n$ and $a$, the LS-SCAD estimator of $\boldsymbol{\beta}$ is

$$\widehat{\boldsymbol{\beta}}_n \equiv \widehat{\boldsymbol{\beta}}(\lambda_n; a) = \arg\min Q_n(\mathbf{b}; \lambda_n, a).$$

We write $\widehat{\boldsymbol{\beta}}_n = (\widehat{\boldsymbol{\beta}}'_{1n}, \widehat{\boldsymbol{\beta}}'_{2n})'$ the way we partition $\boldsymbol{\beta}$ into $\boldsymbol{\beta}_1$ and $\boldsymbol{\beta}_2$.

## 3. Asymptotic properties of the LS-SCAD estimator

In this section we state the results on the asymptotic properties of the LS-SCAD estimator. Results for the case of fixed design are slightly different from those for the case of random design. We state them separately.

For convenience, the main assumptions required for conclusions in this section are listed here. (A0) through (A4) are for fixed covariates. Let $\rho_{n,1}$ be the smallest eigenvalue of $n^{-1}\mathbb{X}'\mathbb{X}$. $\pi_{n,k_n}$ and $\omega_{n,m_n}$ are the largest eigenvalues of $n^{-1}\mathbb{X}'_1\mathbb{X}_1$ and $n^{-1}\mathbb{X}'_2\mathbb{X}_2$, respectively. Let $\mathbf{X}'_{i1} = (X_{i1}, \ldots, X_{ik_n})$ and $\mathbf{X}'_{i2} = (X_{i,k_n+1}, \ldots, X_{ip_n})$.

(A0) (a) $\varepsilon_i$'s are i.i.d with mean 0 and variance $\sigma^2$;
(b) For any $j \in \{1, \ldots, p_n\}$, $\|\mathbb{X}_{\cdot j}\|^2 = n$.
(A1) (a) $\lim_{n\to\infty} \sqrt{k_n}\lambda_n/\sqrt{\rho_{n,1}} = 0$;
(b) $\lim_{n\to\infty} \sqrt{p_n}/\sqrt{n\rho_{n,1}} = 0$.
(A2) (a) $\lim_{n\to\infty} \sqrt{k_n}\lambda_n/(\sqrt{\rho_{n,1}} \min_{1\leq j\leq k_n} |\beta_j|) = 0$;
(b) $\lim_{n\to\infty} \sqrt{p_n}/(\sqrt{n\rho_{n,1}} \min_{1\leq j\leq k_n} |\beta_j|) = 0$;
(c) $\lim_{n\to\infty} \sqrt{p_n/n}/\rho_{n,1} = 0$.
(A3) $\lim_{n\to\infty} \sqrt{\max(\pi_{n,k_n}, \omega_{n,m_n})p_n}/(\sqrt{n}\rho_{n,1}\lambda_n) = 0$.
(A4) $\lim_{n\to\infty} \max_{1\leq i\leq n} \mathbf{X}'_{i1}(\sum_{i=1}^n \mathbf{X}_{i1}\mathbf{X}'_{i1})^{-1}\mathbf{X}_{i1} = 0$.

For random covariates, we require conditions (B0) through (B3). Suppose $(\mathbf{X}'_i, \varepsilon_i)$'s are independent and identically distributed as $(\mathbf{X}', \varepsilon) = (X_1, \ldots, X_{p_n}, \varepsilon)$. Analogous to the fixed design case, $\rho_1$ denotes the smallest eigenvalue of $E[\mathbf{X}\mathbf{X}']$. Also $\pi_{k_n}$ and $\omega_{m_n}$ are the largest eigenvalues of $E[\mathbf{X}_{i1}\mathbf{X}'_{i1}]$ and $E[\mathbf{X}_{i2}\mathbf{X}'_{i2}]$, respectively.

(B0) $(\mathbf{X}'_i, \varepsilon_i) = (X_{i1}, \ldots, X_{ip_n}, \varepsilon_i), i = 1, \ldots, n$ are i.i.d. with
   (a) $E[X_{ij}] = 0$, $\text{Var}(X_{ij}) = 1$;
   (b) $E[\varepsilon|\mathbf{X}] = 0$, $\text{Var}(\varepsilon|\mathbf{X}) = \sigma^2$.
(B1) (a) $\lim_{n\to\infty} p_n^2/(n\rho_1^2) = 0$;
   (b) $\lim_{n\to\infty} k_n \lambda_n^2/\rho_1 = 0$.
(B2) (a) $\lim_{n\to\infty} \sqrt{p_n}/(\sqrt{n\rho_1} \min_{1\le j \le k_n} |\beta_j|) = 0$;
   (b) $\lim_{n\to\infty} \lambda_n \sqrt{k_n}/(\sqrt{\rho_1} \min_{1\le j \le k_n} |\beta_j|) = 0$.
(B3)
$$\lim_{n\to\infty} \frac{\sqrt{\max(\pi_{k_n}, \omega_{m_n})p_n}}{\sqrt{n}\rho_1\lambda_n} = 0.$$

**Theorem 1** (Consistency in the fixed design setting). *Under (A0)–(A1),*
$$\|\widehat{\boldsymbol{\beta}}_n - \boldsymbol{\beta}\| \xrightarrow{P} 0 \quad as\ n \to \infty.$$

A similar result holds for the random design case.

**Theorem 2** (Consistency in the random design setting). *Suppose that there exists an absolute constant $M_4$ such that for all $n$, $\max_{1\le j\le p_n} E[X_j^4] \le M_4 < \infty$. Then under (B0)–(B1),*
$$\|\widehat{\boldsymbol{\beta}}_n - \boldsymbol{\beta}\| \xrightarrow{P} 0 \quad as\ n \to \infty.$$

For consistency, $\lambda_n$ has to be kept small so that the SCAD penalty would not introduce any bias asymptotically. Note that in both design settings, the restriction on the penalty parameter $\lambda_n$ does not involve $m_n$, the number of trivial covariates. This is shared by the $L_q(0 < q < 1)$-penalized estimators in Huang, Horowitz and Ma [4]. However, unlike the bridge estimators, no upper bound requirement is imposed on the components of $\boldsymbol{\beta}_1$, since the derivative of the SCAD penalty vanishes beyond a certain interval while that of the $L_q$ penalty does not. In the fixed design case, (A1.b) is needed for model identifiability, as required in the classical regression. For the random design case, a stricter requirement on $p_n$ is entailed by the need of the convergence of $n^{-1}\mathbb{X}'\mathbb{X}$ to $E[\mathbf{X}\mathbf{X}']$ in the Frobenius norm.

The next two theorems state that the LS-SCAD estimator is consistent for variable selection.

**Theorem 3** (Variable selection in the fixed design setting). *Under (A0)–(A3), $\widehat{\boldsymbol{\beta}}_{2n} = \mathbf{0}_{m_n}$ with probability tending to 1.*

**Theorem 4** (Variable selection in the random design setting). *Suppose there exists an absolute constant $M$ such that $\max_{1\le j\le p_n} |X_j| \le M < \infty$. Then under (B0)–(B3), $\widehat{\boldsymbol{\beta}}_{2n} = \mathbf{0}_{m_n}$ with probability tending to 1.*

(A2.a) and (A2.b) are identical to (A1.a) and (A1.b), respectively, provided that
$$\liminf_{n\to\infty} \min_{1\le j\le k_n} |\beta_j| > 0.$$

(B2) has a requirement for $\min_{1\le j\le k_n} |\beta_j|$ similar to (A2). (A3) concerns the largest eigenvalues of $n^{-1}\mathbb{X}'_1\mathbb{X}_1$ and $n^{-1}\mathbb{X}'_2\mathbb{X}_2$. Due to the standardization of covariates,
$$\pi_{n,k_n} \le k_n \text{ and } \omega_{n,m_n} \le m_n.$$

So (A3) is implied by

$$\lim_{n\to\infty} \frac{p_n}{\sqrt{n}\rho_{n,1}\lambda_n} = 0.$$

Likewise, (B3) can be replaced with

$$\lim_{n\to\infty} \frac{p_n}{\sqrt{n}\rho_1\lambda_n} = 0.$$

Both (A3) and (B3) require $\lambda_n$ not to converge too fast to 0 in order for the estimator to be able to "discover" the trivial covariates. It may be of concern if there are $\lambda_n$'s that simultaneously satisfy (A1)–(A3) (in the random design setting (B1)–(B3)) under certain conditions. When $\liminf \rho_{n,1} > 0$ and $\liminf_{n\to\infty} \min_{1\leq j\leq k_n} |\beta_j| > 0$, it can be checked that there exists $\lambda_n$ that meets both (A2) and (A3) as long as $p_n = o(n^{1/3})$. If we further know either that $k_n$ is fixed, or that the largest eigenvalue of $n^{-1}\mathbb{X}'\mathbb{X}$ is bounded from above, as is assumed in Fan and Peng [2], $p_n = o(n^{1/2})$ is sufficient. When both of these are true, $p_n = o(n)$ is adequate for the existence of such $\lambda_n$'s. Similar conclusions hold for the random design case except that $p_n = o(n^{1/2})$ is indispensable there.

The advantage of the SCAD penalty is that once the trivial covariates have been correctly picked out, regression with or without the SCAD penalty will make no difference to the nontrivial covariates. So it is expected that $\widehat{\boldsymbol{\beta}}_{1n}$ is asymptotically normally distributed. Let $\{\mathbf{A}_n, n = 1, 2, \ldots\}$ be a sequence of matrices of dimension $d \times k_n$ with full row rank.

**Theorem 5** (Asymptotic normality in the fixed design setting). *Under (A0)–(A4),*

$$\sqrt{n}\,\Sigma_n^{-1/2}\mathbf{A}_n(\widehat{\boldsymbol{\beta}}_{1n} - \boldsymbol{\beta}_1) \xrightarrow{D} N(\mathbf{0}_d, I_d),$$

where $\Sigma_n = \sigma^2 \mathbf{A}_n (\sum_{i=1}^n \mathbf{X}_{i1}\mathbf{X}'_{i1}/n)^{-1}\mathbf{A}'_n$.

**Theorem 6** (Asymptotic normality in the random design setting). *Suppose that there exists an absolute constant $M$ such that $\max_{1\leq j\leq p_n} \|X_j\| \leq M < \infty$ and a $\sigma_4$ such that $E[\varepsilon^4|\mathbb{X}_{11}] \leq \sigma_4 < \infty$ for all $n$. Then under (B0)–(B3),*

$$n^{-1/2}\Sigma_n^{-1/2}\mathbf{A}_n E^{-1/2}[\mathbf{X}_{i1}\mathbf{X}'_{i1}] \sum_{i=1}^n \mathbf{X}_{i1}\mathbf{X}'_{i1}(\widehat{\boldsymbol{\beta}}_{1n} - \boldsymbol{\beta}_1) \xrightarrow{D} N(\mathbf{0}_d, I_d),$$

where $\Sigma_n = \sigma^2 \mathbf{A}_n \mathbf{A}'_n$.

For the random design the assumptions for asymptotic normality are no more than those for variable selection. While for the fixed design, a Lindeberg-Feller condition (A4) is needed in addition to (A0)–(A3).

## 4. Computation

We use the algorithm of Hunter and Li [6] to compute the LS-SCAD estimator for a given $\lambda_n$ and $a$. This algorithm approximates a nonconvex target function with a convex function locally at each iteration step. We also describe the steps to compute the approximate standard error of the estimator.

### 4.1. Computation of the LS-SCAD estimator

Given $\lambda_n$ and $a$ the target function to be minimized is

$$Q_n(\mathbf{b}; \lambda_n, a) = \sum_{i=1}^{n}(Y_i - \mathbf{X}_i'\mathbf{b})^2 + n\sum_{j=1}^{p_n} p_{\lambda_n}(b_j; a).$$

Hunter and Li [6] proposes to minimize its approximation

$$Q_{n,\xi}(\mathbf{b}; \lambda_n, a) = \sum_{i=1}^{n}(Y_i - \mathbf{X}_i'\mathbf{b})^2 + n\sum_{j=1}^{p_n} p_{\lambda_n,\xi}(b_j; a)$$

$$= \sum_{i=1}^{n}(Y_i - \mathbf{X}_i'\mathbf{b})^2 + n\sum_{j=1}^{p_n}\left(p_{\lambda_n}(b_j;a) - \xi\int_0^{|b_j|}\frac{p'_{\lambda_n}(t;a)}{\xi+t}dt\right)$$

Around $\mathbf{b}_{(k)} = (b_{(k),1}, \ldots, b_{(k),p_n})'$, it can be approximated by

$$S_{k,\xi}(\mathbf{b}; \lambda_n, a) = \sum_{i=1}^{n}(Y_i - \mathbf{X}_i'\mathbf{b})^2$$

$$+ n\sum_{j=1}^{p_n}\left[p_{\lambda_n,\xi}(b_{(k),j};a) + \frac{p'_{\lambda_n}(|b_{(k),j}|+;a)}{2(\xi + |b_{(k),j}|)}(b_j^2 - b_{(k),j}^2)\right],$$

where $\xi$ is a very small perturbation to prevent any component of the estimate from getting stuck at 0. Therefore the one-step estimator starting from $\mathbf{b}_{(k)}$ is

$$\mathbf{b}_{(k+1)} = (\mathbb{X}'\mathbb{X} + n\mathbf{D}_\xi(\mathbf{b}_{(k)}; \lambda_n, a))^{-1}\mathbb{X}'\mathbf{Y},$$

where $\mathbf{D}_\xi(\mathbf{b}_{(k)}; \lambda_n, a)$ is the diagonal matrix whose diagonal elements are $\frac{1}{2}p'_{\lambda_n} \times (|b_{(k),j}|+;a)/(\xi + |b_{(k),j}|), j = 1, \ldots, p_n$. Given the tolerance $\tau$, convergence is claimed when

$$\left|\frac{\partial Q_{n,\xi}(\mathbf{b})}{\partial b_j}\right| < \frac{\tau}{2}, \quad \forall j = 1, \ldots, p_n.$$

And finally the $b_j$'s that satisfy

$$\left|\frac{\partial Q_{n,\xi}(\mathbf{b})}{\partial b_j} - \frac{\partial Q_n(\mathbf{b})}{\partial b_j}\right| = \frac{n\xi p'_{\lambda_n}(|b_j|+;a)}{\xi + |b_j|} > \frac{\tau}{2}$$

are set to 0. A good starting point would be $\mathbf{b}_{(0)} = \widehat{\boldsymbol{\beta}}_{\text{LS}}$, the least squares estimator.

The perturbation $\xi$ should be kept small so that difference between $Q_{n,\xi}(\cdot)$ and $Q_n(\cdot)$ is negligible. Hunter and Li [6] suggests using

$$\xi = \frac{\tau}{2n\lambda_n}\min\{|b_{(0),j}| : b_{(0),j} \neq 0\}.$$

### 4.2. Standard errors

The standard errors for the nonzero coefficient estimates can be obtained via the approximation

$$\frac{\partial S_\xi(\widehat{\boldsymbol{\beta}}_{1n}; \lambda, a)}{\partial \widehat{\boldsymbol{\beta}}_{1n}} \approx \frac{\partial S_\xi(\boldsymbol{\beta}_1; \lambda_n, a)}{\partial \boldsymbol{\beta}_1} + \frac{\partial^2 S_\xi(\boldsymbol{\beta}_1; \lambda_n, a)}{\partial \boldsymbol{\beta}_1 \partial \boldsymbol{\beta}_1'}\left(\widehat{\boldsymbol{\beta}}_{1n} - \boldsymbol{\beta}_1\right).$$

So

$$\widehat{\boldsymbol{\beta}}_{1n} - \boldsymbol{\beta}_1 \approx - \left( \frac{\partial^2 S_\xi(\boldsymbol{\beta}_1; \lambda_n, a)}{\partial \boldsymbol{\beta}_1 \partial \boldsymbol{\beta}_1'} \right)^{-1} \frac{\partial S_\xi(\boldsymbol{\beta}_1; \lambda_n, a)}{\partial \boldsymbol{\beta}_1}$$

$$\approx - \left( \frac{\partial^2 S_\xi(\widehat{\boldsymbol{\beta}}_{1n}; \lambda_n, a)}{\partial \widehat{\boldsymbol{\beta}}_{1n} \partial \widehat{\boldsymbol{\beta}}_{1n}'} \right)^{-1} \frac{\partial S_\xi(\widehat{\boldsymbol{\beta}}_{1n}; \lambda_n, a)}{\partial \widehat{\boldsymbol{\beta}}_{1n}}.$$

Since

$$\frac{\partial S_\xi(\widehat{\boldsymbol{\beta}}_{1n}; \lambda_n, a)}{\partial \widehat{\beta}_j} = -2\mathbb{X}_{\cdot j}' \mathbf{Y} + 2\mathbb{X}_{\cdot j}' \mathbb{X}_1 \widehat{\boldsymbol{\beta}}_{1n} + n\frac{\widehat{\beta}_j p_{\lambda_n}'(|\widehat{\beta}_j|; a)}{\xi + |\widehat{\beta}_j|}$$

$$= \sum_{i=1}^n \left[ -2 X_{ij} Y_i + 2 X_{ij} \mathbf{X}_{i1}' \widehat{\boldsymbol{\beta}}_{1n} + \frac{\widehat{\beta}_j p_{\lambda_n}'(|\widehat{\beta}_j|; a)}{\xi + |\widehat{\beta}_j|} \right],$$

$$\triangleq 2 \sum_{i=1}^n U_{ij}(\xi; \lambda_n, a),$$

letting $U_{ij} = U_{ij}(\xi; \lambda_n, a)$, we have, for $j, l = 1, \ldots, k_n$,

$$\text{Cov} \left( n^{-1/2} \frac{\partial S_\xi(\widehat{\boldsymbol{\beta}}_{1n}; \lambda_n, a)}{\partial \widehat{\beta}_j}, n^{-1/2} \frac{\partial S_\xi(\widehat{\boldsymbol{\beta}}_{1n}; \lambda_n, a)}{\partial \widehat{\beta}_l} \right)$$

$$\approx \frac{4}{n} \sum_{i=1}^n U_{ij} U_{il} - \frac{4}{n^2} \sum_{i=1}^n U_{ij} \sum_{i=1}^n U_{il}.$$

Let $\mathbb{C} = (C_{jl}, j, l = 1, \ldots, k_n)$, where

$$C_{jl} = \frac{1}{n} \sum_{i=1}^n U_{ij} U_{il} - \frac{1}{n^2} \sum_{i=1}^n U_{ij} \sum_{i=1}^n U_{il}.$$

The variance-covariance matrix of the estimates can be approximated by

$$\widehat{\text{Cov}(\widehat{\boldsymbol{\beta}}_{1n})} \equiv n(\mathbb{X}_1' \mathbb{X}_1 + n\mathbf{D}_\xi(\widehat{\boldsymbol{\beta}}_{1n}; \lambda_n, a))^{-1} \mathbb{C} \, (\mathbb{X}_1' \mathbb{X}_1 + n\mathbf{D}_\xi(\widehat{\boldsymbol{\beta}}_{1n}; \lambda_n, a))^{-1}.$$

### 4.3. Selection of $\lambda_n$

The above computational algorithm is for the case when $\lambda_n$ and $a$ are specified. In data analysis, they can be selected by minimizing the generalized cross validation score, which is defined to be

$$\text{GCV}(\lambda_n, a) = \frac{\|\mathbf{Y} - \mathbb{X}_1 \widehat{\boldsymbol{\beta}}_{1n}\|^2/n}{(1 - p(\lambda_n, a)/n)^2},$$

where

$$p(\lambda_n, a) = \text{tr} \left[ \mathbb{X}_1 \left( \mathbb{X}_1' \mathbb{X}_1 + n\mathbf{D}_0(\widehat{\boldsymbol{\beta}}_{1n}; \lambda_n, a) \right)^{-1} \mathbb{X}_1' \right]$$

is the number of effective parameters and $\mathbf{D}_0(\widehat{\boldsymbol{\beta}}_{1n}; \lambda_n, a)$ is a submatrix of the diagonal matrix $\mathbf{D}_\xi(\widehat{\boldsymbol{\beta}}_n; \lambda_n, a)$ with $\xi = 0$. By submatrix, we mean the diagonal of $\mathbf{D}_0(\widehat{\boldsymbol{\beta}}_{1n}; \lambda_n, a)$ only contains the elements corresponding to the nontrivial components in $\widehat{\boldsymbol{\beta}}$. Note that here $\mathbb{X}_1$ also only includes the columns of which the corresponding elements of $\widehat{\boldsymbol{\beta}}_n$ are non-vanishing.

The requirement that $a > 2$ is implied by the SCAD penalty function. Simulation suggests that the generalized cross validation score does not change much with $a$ given $\lambda$. So to improve computing efficiency, we fix $a = 3.7$, as suggested by Fan and Li [1].

## 5. Simulation studies

In this section we illustrate the LS-SCAD estimator's finite sample properties with a simulated example.

We simulate covariates $\mathbf{X}_i, i = 1, \ldots, n$ from the multivariate normal distributions with mean 0 and

$$\mathrm{Cov}(X_{ij}, X_{il}) = \rho^{|j-l|}, 1 \le j, l \le p,$$

The response $Y_i$ is computed as

$$Y_i = \sum_{j=1}^{p} X_{ij}\beta_j + \varepsilon_i, \qquad i = 1, \ldots, n.$$

where $\beta_j = j, 1 \le j \le 4$, $\beta_j = 0, 5 \le j \le p$, and $\varepsilon_i$'s are sampled from $N(0,1)$. For each $(n, p, \rho) \in \{(100, 10), (500, 40)\} \times \{0, 0.2, 0.5, 0.8\}$, we generated $N = 400$ data sets and use the algorithm in Section 4 to compute the LS-SCAD estimator. We set the tolerance $\tau$ described in Section 4.1 at $10^{-5}$. For comparison we also apply the ordinary least square (LS) method, the ordinary least square method with model selection based on AIC (abbreviated as AIC), and the ordinary least squares assuming that $\beta_j = 0$ for $j \ge 5$ are known beforehand (ORA). Note that this last estimator (ORA) is not feasible in a real data analysis setting. We use it here as a benchmark in the comparisons.

The results are summarized in Tables 1 and 2. Columns 4 through 7 in Table 1 are the biases of the estimates of $\beta_j, j = 1, \ldots, 4$ respectively. In the parentheses following each of them are the standard deviations of these estimates. Column 8 ($\overline{K}$) lists the numbers of estimates of $\beta_j, 5 \le j \le p$ that are 0, averaged over 400 replications, and their modes are given in Column 9 ($\widetilde{K}$). For LS, an estimate is set to be 0 if it lies within $[-10^{-5}, 10^{-5}]$.

In Table 1, we see that the LS-SCAD estimates of the nontrivial coefficients have biases and standard errors comparable to the ORA estimates. This is in line with Theorems 5 and 6. The average numbers of nonzero estimates for $\beta_j(j > 4)$, $\overline{K}$, with respect to LS-SCAD are close to $p$, the true number of nonzero coefficients among $\beta_j(j > 4)$. As the true number of trivial covariates increases, the LS-SCAD estimator may be able to discover more trivial ones than AIC. However, there is more variability in the number of trivial covariates discovered via LS-SCAD than that via AIC.

Table 2 gives the averages of the estimated standard errors of $\widehat{\beta}_j, 1 \le j \le 4$ using the SCAD method over the 400 replications. They are obtained based on the approach described in Section 4.2. They are slightly smaller than the sampling standard deviations of $\widehat{\beta}_j, 1 \le j \le 4$, which are given in parentheses in the rows for LS-SCAD.

Suppose for a data set the estimate of $\boldsymbol{\beta}$ via one of these four approaches is $\widehat{\boldsymbol{\beta}}$, then the average model error (AME) regarding this approach is computed as $n^{-1} \sum_{i=1}^{n} [\mathbf{X}'_i (\widehat{\boldsymbol{\beta}}_n - \boldsymbol{\beta})]^2$. Box plots for these AME's are given in Figure 1.

TABLE 1
*Simulation example 1, comparison of estimators*

| $(n,p)$ | $\rho$ | Estimator | $\beta_1$ | $\beta_2$ | $\beta_3$ | $\beta_4$ | $\overline{K}$ | $\widetilde{K}$ |
|---|---|---|---|---|---|---|---|---|
| (100, 10) | 0 | LS | .0007 (.1112) | −.0034 (.0979) | −.0064 (.1127) | −.0024 (.1091) | 0 | 0 |
| | | ORA | .0008 (.1074) | −.0054 (.0936) | −.0057 (.1072) | −.0007 (.1040) | 6 | 6 |
| | | AIC | .0007 (.1083) | −.0026 (.1033) | −.0060 (.1156) | −.0019 (.1181) | 4.91 | 5 |
| | | SCAD | −.0006 (.1094) | −.0037 (.0950) | −.0058 (.1094) | −.0014 (.1060) | 4.62 | 5 |
| | 0.2 | LS | −.0003 (.1051) | −.0028 (.1068) | .0093 (.1157) | .0037 (.1103) | 0 | 0 |
| | | ORA | −.0005 (.1010) | −.0031 (.1035) | .0107 (.1131) | .0020 (.1035) | 6 | 6 |
| | | AIC | −.0002 (.1031) | −.0024 (.1063) | .0107 (.1150) | .0021 (.1079) | 4.95 | 5 |
| | | SCAD | −.0025 (.1035) | −.0026 (.1046) | .0104 (.1141) | .0024 (.1066) | 4.64 | 5 |
| | 0.5 | LS | .0000 (.1177) | −.0007 (.1353) | .0010 (.1438) | .0006 (.1360) | 0 | 0 |
| | | ORA | −.0002 (.1129) | −.0072 (.1317) | .0115 (.1393) | .0022 (.1171) | 6 | 6 |
| | | AIC | −.0003 (.1162) | −.0064 (.1338) | .0114 (.1413) | .0017 (.1294) | 4.91 | 5 |
| | | SCAD | .0035 (.1115) | −.0219 (.1404) | .0135 (.1481) | .0006 (.1293) | 4.78 | 5 |
| | 0.8 | LS | −.0005 (.1916) | −.0229 (.2293) | .0059 (.2319) | .0060 (.2200) | 0 | 0 |
| | | ORA | −.0039 (.1835) | −.0196 (.2197) | .0070 (.2250) | .0092 (.1787) | 6 | 6 |
| | | AIC | −.0021 (.1857) | −.0209 (.2235) | .0063 (.2289) | .0013 (.2072) | 4.85 | 6 |
| | | SCAD | −.0038 (.1868) | −.0197 (.2249) | .0062 (.2280) | .0032 (.2024) | 4.87 | 6 |
| (500, 40) | 0 | LS | .0021 (.0466) | −.0000 (.0475) | −.0010 (.0466) | .0014 (.0439) | 0 | 0 |
| | | ORA | .0027 (.0446) | −.0005 (.0453) | −.0003 (.0448) | .0011 (.0426) | 36 | 36 |
| | | AIC | .0023 (.0460) | −.0003 (.0465) | −.0004 (.0453) | .0016 (.0433) | 29.91 | 30 |
| | | SCAD | .0027 (.0447) | −.0004 (.0454) | −.0004 (.0450) | .0013 (.0429) | 32.22 | 35 |
| | 0.2 | LS | .0018 (.0478) | .0003 (.0478) | −.0014 (.0487) | .0005 (.0437) | 0 | 0 |
| | | ORA | .0003 (.0522) | −.0000 (.0465) | −.0010 (.0517) | .0009 (.0458) | 36 | 36 |
| | | AIC | .0024 (.0473) | .0002 (.0471) | −.0014 (.0475) | .0018 (.0436) | 29.87 | 30 |
| | | SCAD | .0028 (.0461) | .0002 (.0460) | −.0011 (.0475) | .0006 (.0433) | 32.20 | 35 |
| | 0.5 | LS | .0024 (.0542) | .0001 (.0617) | .0050 (.0608) | −.0048 (.0563) | 0 | 0 |
| | | ORA | .0027 (.0526) | .0017 (.0581) | .0033 (.0597) | −.0030 (.0488) | 36 | 36 |
| | | AIC | .0031 (.0537) | .0007 (.0603) | .0037 (.0605) | −.0038 (.0526) | 29.87 | 32 |
| | | SCAD | .0025 (.0528) | .0017 (.0587) | .0034 (.0601) | −.0037 (.0494) | 31.855 | 35 |
| | 0.8 | LS | .0014 (.0788) | −.0012 (.1014) | .0090 (.1000) | −.0077 (.0943) | 0 | 0 |
| | | ORA | .0010 (.0761) | .0017 (.0954) | .0060 (.0983) | −.0044 (.0704) | 36 | 36 |
| | | AIC | .0020 (.0776) | .0003 (.0996) | .0066 (.0995) | −.0071 (.0862) | 29.56 | 30 |
| | | SCAD | .0014 (.0773) | .0018 (.0982) | .0059 (.0990) | −.0050 (.0790) | 29.38 | 35 |

TABLE 2
*Simulated example, standard error estimate*

| $(n,p)$ | (100, 10) | | | | (500, 40) | | | |
|---|---|---|---|---|---|---|---|---|
| $\rho$ | 0 | 0.2 | 0.5 | 0.8 | 0 | 0.2 | 0.5 | 0.8 |
| $se(\widehat{\beta}_1)$ | .0983 | .1005 | .1139 | .1624 | .0442 | .0444 | .0512 | .0735 |
| $se(\widehat{\beta}_2)$ | .0980 | .1028 | .1276 | .2080 | .0443 | .0447 | .0571 | .0940 |
| $se(\widehat{\beta}_3)$ | .0996 | .1027 | .1278 | .2086 | .0442 | .0445 | .0573 | .0940 |
| $se(\widehat{\beta}_4)$ | .0988 | .1006 | .1150 | .1727 | .0441 | .0444 | .0512 | .0764 |

The LS estimator definitely has the worst performance in terms of AME. This becomes more obvious as the number of trivial predictors increases. LS-SCAD outperforms AIC in this respect and is comparable to ORA. But it is also seen that the AME's of LS-SCAD tend to be more diffuse as $\rho$ increases. This is also the result of more spread-out estimates of the number of trivial covariates.

## 6. Concluding remarks

In this paper, we have studied the asymptotic properties of the LS-SCAD estimator when the number of covariates and regression coefficients increases to infinity as

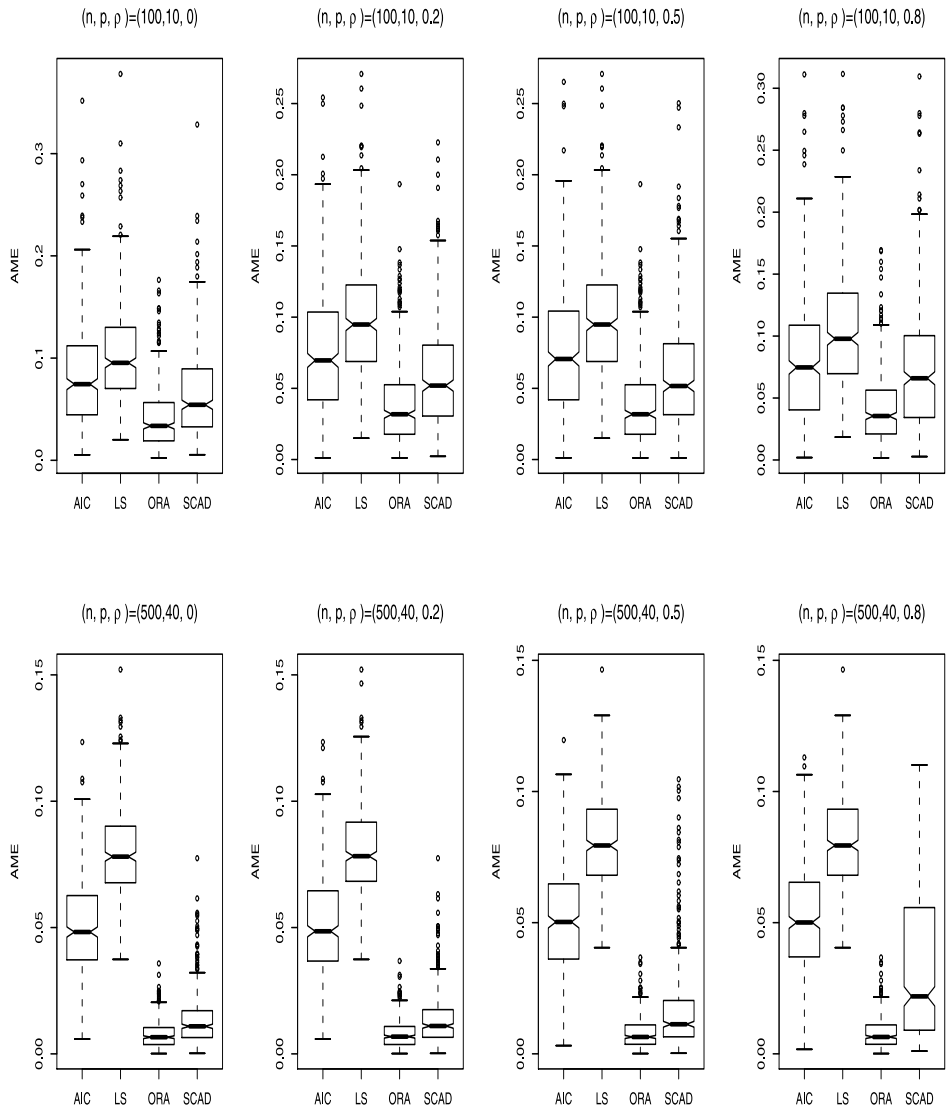

FIG 1. *Box plots of the average model errors for four estimators: AIC, LS, ORA, and LS-SCAD. In the top four panels, $(n, p, \rho) = (100, 10, 0), (100, 10, 0.2), (100, 10, 0.5), (100, 10, 0.8)$; and in the bottom four panels, $(n, p, \rho) = (500, 40, 0), (500, 40, 0.2), (500, 40, 0.5), (500, 40, 0.8)$, where $n$ is the sample size, $p$ is the number of covariates, and $\rho$ is the correlation coefficient used in generating the covariate values.*

$n \to \infty$. We have shown that this estimator can correctly identify zero coefficients with probability converging to one and that the estimators of nonzero coefficients are asymptotically normal and oracle efficient. Our results were obtained under the assumption that the number of parameters is smaller than the sample size. They are not applicable when the number of parameters is greater than the sample size, which arises in microarray gene expression studies. In general, the condition that $p < n$ is needed for identification of the regression parameter and consistent variable selection. To achieve consistent variable selection in the "large $p$, small $n$" case, certain conditions are required for the design matrix. For example, Huang et al. [4]

showed that, under a partial orthogonality assumption in which the covariates of the zero coefficients are uncorrelated or only weakly correlated with the covariates of nonzero coefficients, then the univariate bridge estimators are consistent for variable selection under appropriate conditions. This result also holds for the univariate LS-SCAD estimator. Indeed, under the partial orthogonality condition, it can be shown that the simple univariate regression estimator can be used to consistently distinguish between nonzero and zero coefficients. Finally, we note that our results are only valid for a fixed sequence of penalty parameters $\lambda_n$. It is an interesting and difficult problem to show that the asymptotic oracle property also holds for $\lambda_n$ determined by cross validation.

**Appendix**

We now give the proofs of the results stated in Section 3.

*Proof of Theorem 1.* By the definition of $\widehat{\boldsymbol{\beta}}_n$, it is necessary that $Q_n(\widehat{\boldsymbol{\beta}}_n) \leq Q_n(\boldsymbol{\beta})$. It follows that

$$
\begin{aligned}
0 &\geq \left\|\mathbb{X}(\widehat{\boldsymbol{\beta}}_n - \boldsymbol{\beta})\right\|^2 - 2\boldsymbol{\varepsilon}'\mathbb{X}(\widehat{\boldsymbol{\beta}}_n - \boldsymbol{\beta}) + n\sum_{j=1}^{k_n}\left[p_{\lambda_n}(\widehat{\beta}_j;a) - p_{\lambda_n}(\beta_j;a)\right] \\
&\geq \left\|\mathbb{X}(\widehat{\boldsymbol{\beta}}_n - \boldsymbol{\beta})\right\|^2 - 2\boldsymbol{\varepsilon}'\mathbb{X}(\widehat{\boldsymbol{\beta}}_n - \boldsymbol{\beta}) - 2^{-1}n(a+1)k_n\lambda_n^2 \\
&= \left\|[\mathbb{X}'\mathbb{X}]^{1/2}(\widehat{\boldsymbol{\beta}}_n - \boldsymbol{\beta}) - [\mathbb{X}'\mathbb{X}]^{-1/2}\mathbb{X}'\boldsymbol{\varepsilon}\right\|^2 \\
&\quad - \boldsymbol{\varepsilon}'\mathbb{X}[\mathbb{X}'\mathbb{X}]^{-1}\mathbb{X}'\boldsymbol{\varepsilon} - 2^{-1}n(a+1)k_n\lambda_n^2.
\end{aligned}
$$

By the $C_r$-inequality (Loéve [8], page 155),

$$
\begin{aligned}
\left\|[\mathbb{X}'\mathbb{X}]^{1/2}(\widehat{\boldsymbol{\beta}}_n - \boldsymbol{\beta})\right\|^2 &\leq 2\left\|[\mathbb{X}'\mathbb{X}]^{1/2}(\widehat{\boldsymbol{\beta}}_n - \boldsymbol{\beta}) - [\mathbb{X}'\mathbb{X}]^{-1/2}\mathbb{X}'\boldsymbol{\varepsilon}\right\|^2 + 2\boldsymbol{\varepsilon}'\mathbb{X}[\mathbb{X}'\mathbb{X}]^{-1}\mathbb{X}'\boldsymbol{\varepsilon} \\
&\leq 4\boldsymbol{\varepsilon}'\mathbb{X}[\mathbb{X}'\mathbb{X}]^{-1}\mathbb{X}'\boldsymbol{\varepsilon} + n(a+1)k_n\lambda_n^2.
\end{aligned}
$$

In the fixed design,

$$
\begin{aligned}
\boldsymbol{\varepsilon}'\mathbb{X}[\mathbb{X}^{(n)'}\mathbb{X}]^{-1}\mathbb{X}'\boldsymbol{\varepsilon} &= E\left[\boldsymbol{\varepsilon}'\mathbb{X}[\mathbb{X}^{(n)'}\mathbb{X}]^{-1}\mathbb{X}'\boldsymbol{\varepsilon}\right]O_P(1) \\
&= \sigma^2 \mathrm{tr}(\mathbb{X}[\mathbb{X}'\mathbb{X}]^{-1}\mathbb{X}')O_P(1) \\
&= p_n O_P(1).
\end{aligned}
$$

Since

$$\left\|[\mathbb{X}'\mathbb{X}]^{1/2}(\widehat{\boldsymbol{\beta}}_n - \boldsymbol{\beta})\right\|^2 \geq n\rho_{n,1}\|\widehat{\boldsymbol{\beta}}_n - \boldsymbol{\beta}\|^2,$$

we have

$$\|\widehat{\boldsymbol{\beta}}_n - \boldsymbol{\beta}\| = O_P\left(\frac{\sqrt{p_n}}{\sqrt{n\rho_{n,1}}} + \frac{\sqrt{k_n}\lambda_n}{\sqrt{\rho_{n,1}}}\right) = o_P(1). \qquad \square$$

*Proof of Theorem 2.* Let $\mathbf{A}^{(n)} = (A_{jk}^{(n)})_{j,k=1,\ldots,p_n}$ with $A_{jk}^{(n)} = n^{-1}\sum_{i=1}^{n}X_{ij}X_{ik} - E[X_{ij}X_{ik}]$. Let $\rho_1(\mathbf{A}^{(n)})$ and $\rho_{p_n}(\mathbf{A}^{(n)})$ be the smallest and largest of the eigenvalues of $\mathbf{A}^{(n)}$, respectively. Then by Theorem 4.1 in Wang and Jia [13],

$$\rho_1(\mathbf{A}^{(n)}) \leq \rho_{n,1} - \rho_1 \leq \rho_{p_n}(\mathbf{A}).$$

By the Cauchy inequality and the properties of eigenvalues of symmetric matrices,

$$\max(|\rho_1(\mathbf{A}^{(n)})|, |\rho_{p_n}(\mathbf{A}^{(n)})|) \leq \|\mathbf{A}^{(n)}\|.$$

When (B1.a) holds, $\|\mathbf{A}^{(n)}\| = o_P(\rho_1) = o_P(1)$, as is seen for any $\xi > 0$,

$$P(\|\mathbf{A}^{(n)}\|^2 \geq \xi \rho_1{}^2) \leq \frac{E\|\mathbf{A}^{(n)}\|^2}{\xi \rho_1{}^2} \leq \frac{p_n^2}{\xi \rho_1{}^2} \sup_{1 \leq j, k \leq p_n} \mathrm{Var}(A_{jk}^{(n)}) \leq \frac{p_n^2}{n\xi \rho_1{}^2} M_4.$$

Since $\rho_1 > 0$ holds for all $n$, $n^{-1}\mathbb{X}'\mathbb{X}$ is invertible with probability tending to 1.

Following the argument for the fixed design case, with probability tending to 1,

$$\left\|[\mathbb{X}'\mathbb{X}]^{1/2}(\widehat{\boldsymbol{\beta}}_n - \boldsymbol{\beta})\right\|^2 \leq 4\boldsymbol{\varepsilon}'\mathbb{X}[\mathbb{X}'\mathbb{X}]^{-1}\mathbb{X}'\boldsymbol{\varepsilon} - n(a+1)k_n\lambda_n^2.$$

In the random design setting,

$$E\left[\boldsymbol{\varepsilon}'\mathbb{X}[\mathbb{X}'\mathbb{X}]^{-1}\mathbb{X}'\boldsymbol{\varepsilon} \,\Big|\, \|\mathbf{A}^{(n)}\|^2 < \frac{1}{2}\rho_1{}^2\right]$$
$$= \sigma^2 E\left[\mathrm{tr}(\mathbb{X}[\mathbb{X}'\mathbb{X}]^{-1}\mathbb{X}') \,\Big|\, \|\mathbf{A}^{(n)}\|^2 < \frac{1}{2}\rho_1{}^2\right]$$
$$= \sigma^2 p_n.$$

The rest of the argument remains the same as for the fixed design case and leads to

$$\|\widehat{\boldsymbol{\beta}}_n - \boldsymbol{\beta}\| = O_P\left(\frac{\sqrt{p_n}}{\sqrt{n\rho_1}} + \frac{\sqrt{k_n}\lambda_n}{\sqrt{\rho_1}}\right) = o_P(1). \qquad \square$$

**Lemma 1** (Convergency rate in the fixed design setting). *Under (A0)–(A2), $\|\widehat{\boldsymbol{\beta}}_n - \boldsymbol{\beta}\| = O_P(\sqrt{p_n/n}/\rho_{n,1})$.*

*Proof.* In the proof of consistency, we have

$$\|\widehat{\boldsymbol{\beta}}_n - \boldsymbol{\beta}\| = O_P(u_n), \quad \text{where } u_n = \lambda_n\sqrt{k_n/\rho_{n,1}} + \sqrt{p_n/(n\rho_{n,1})}.$$

For any $L_1$, provided that $\|\mathbf{b} - \boldsymbol{\beta}\| \leq 2^{L_1}u_n$,

$$\min_{1 \leq j \leq k_n} |b_j| \geq \min_{1 \leq j \leq k_n} |\beta_j| - 2^{L_1}u_n.$$

If (A2) holds, then for $n$ sufficiently large, $u_n/\min_{1 \leq j \leq k_n} |\beta_j| < 2^{-L_1-1}$. It follows that

$$\min_{1 \leq j \leq k_n} |b_j| \geq \min_{1 \leq j \leq k_n} |\beta_j|/2,$$

which further implies than $\min_{1 \leq j \leq k_n} |b_j| > a\lambda_n$ for n sufficiently large (assume $\liminf_{n \to \infty} k_n > 0$).

Let $\{h_n\}$ be a sequence converging to 0. As in the proof of of Theorem 3.2.5 of Van der Vaart and Wellner [12], decompose $\mathcal{R}^{p_n} \backslash \{\mathbf{0}_{p_n}\}$ into shells $\{S_{n,l}, l \in \bar{\mathcal{Z}}\}$ where $S_{n,l} = \{\mathbf{b} : 2^{l-1}h_n \leq \|\mathbf{b} - \boldsymbol{\beta}\| < 2^l h_n\}$. For $\mathbf{b} \in S_{n,l}$ such that $2^l h_n \leq 2^{L_1}u_n$,

$$Q_n(\mathbf{b}) - Q_n(\boldsymbol{\beta}) = (\mathbf{b} - \boldsymbol{\beta})'\mathbb{X}'\mathbb{X}(\mathbf{b} - \boldsymbol{\beta}) - 2\boldsymbol{\varepsilon}'\mathbb{X}(\mathbf{b} - \boldsymbol{\beta})$$
$$+ n\sum_{j=1}^{p_n} p_{\lambda_n}(b_j; a) - n\sum_{j=1}^{p_n} p_{\lambda_n}(\beta_j; a)$$
$$= (\mathbf{b} - \boldsymbol{\beta})'\mathbb{X}'\mathbb{X}(\mathbf{b} - \boldsymbol{\beta}) - 2\boldsymbol{\varepsilon}'\mathbb{X}(\mathbf{b} - \boldsymbol{\beta})$$
$$\triangleq I_{n1} + I_{n2},$$

and

$$I_{n1} \geq n\rho_{n,1}\|\mathbf{b} - \boldsymbol{\beta}\|^2 \geq 2^{2(l-1)}h_n^2 n\rho_{n,1}.$$

Thus

$$P\left(\|\widehat{\boldsymbol{\beta}}_n - \boldsymbol{\beta}\| \geq 2^L h_n\right)$$
$$\leq o(1) + \sum_{\substack{l > L, \\ 2^l h_n \leq 2^{L_1} u_n}} P\left(\widehat{\boldsymbol{\beta}}_n \in S_{n,l}\right)$$
$$\leq o(1) + \sum_{\substack{l > L, \\ 2^l h_n \leq 2^{L_1} u_n}} P\left(\inf_{\mathbf{b} \in S_{n,l}} Q_n(\mathbf{b}) \leq Q_n(\boldsymbol{\beta})\right)$$
$$\leq o(1) + \sum_{\substack{l > L, \\ 2^{l-1} h_n \leq 2^{L_1} u_n}} P\left(\sup_{\mathbf{b} \in S_{n,l}} \boldsymbol{\varepsilon}' \mathbb{X}(\mathbf{b} - \boldsymbol{\beta}) \geq 2^{2l-3} h_n^2 n \rho_{n,1}\right)$$
$$\leq o(1) + \sum_{\substack{l > L, \\ 2^{l-1} h_n \leq 2^{L_1} u_n}} \frac{E|\sup_{\mathbf{b} \in S_{n,l}} \boldsymbol{\varepsilon}' \mathbb{X}(\mathbf{b} - \boldsymbol{\beta})|}{2^{2l-3} h_n^2 n \rho_{n,1}}$$
$$\leq o(1) + \sum_{l > L} \frac{2^l h_n E^{1/2}[\|\boldsymbol{\varepsilon}' \mathbb{X}\|^2]}{2^{2l-3} h_n^2 n \rho_{n,1}}$$
$$\leq o(1) + \sum_{l > L} \frac{2^l \sqrt{n \sigma^2 p_n}}{2^{2l-3} h_n n \rho_{n,1}},$$

from which we see $\|\widehat{\boldsymbol{\beta}}_n - \boldsymbol{\beta}\| = O_P(\sqrt{p_n/n}/\rho_{n,1})$. □

**Lemma 2** (Convergence rate in the random design setting). *Under (B0)–(B2), $\|\widehat{\boldsymbol{\beta}}_n - \boldsymbol{\beta}\| = O_P(\sqrt{p_n/n}/\rho_1)$.*

*Proof.* Deduction is similar to that of Lemma 1. However, since $\mathbb{X}$ is a random matrix in this case, extra details are needed in the following part. Let $\mathbf{A}^{(n)} = (A_{jk}^{(n)})_{j,k=1,\ldots,p_n}$ with $A_{jk}^{(n)} = \frac{1}{n} \sum_{i=1}^n X_{ij} X_{ik} - E[X_j X_k]$. We have

$$P\left(\|\widehat{\boldsymbol{\beta}}_n - \boldsymbol{\beta}\| \geq 2^L h_n\right)$$
$$\leq \sum_{\substack{l > L \\ 2^l h_n \leq 2^{L_1} u_n}} P\left(\widehat{\boldsymbol{\beta}}_n \in S_{n,l}, \|\mathbf{A}^{(n)}\| \leq \rho_1/2\right) + o(1)$$
$$\leq \sum_{\substack{l > L \\ 2^l h_n \leq 2^{L_1} u_n}} P\left(\inf_{\mathbf{b} \in S_{n,l}} Q_n(\mathbf{b}) \leq Q_n(\boldsymbol{\beta}), \|\mathbf{A}^{(n)}\| \leq \rho_1/2\right) + o(1)$$
$$\leq \sum_{l > L} \frac{2^l h_n E^{1/2}\left[\|\boldsymbol{\varepsilon}' \mathbb{X}\|^2 \big| \|\mathbf{A}\| \leq \rho_1/2\right]}{2^{2l-4} h_n^2 n \rho_1} + o(1).$$

The first inequality follows from (B1.a). This leads to $\|\widehat{\boldsymbol{\beta}}_n - \boldsymbol{\beta}\| = O_P(\sqrt{p_n/n}/\rho_1)$. □

*Proof of Theorem 3.* By Lemma 1, $\|\widehat{\boldsymbol{\beta}}_n - \boldsymbol{\beta}\| \leq \lambda_n$ with probability tending to 1 under (A3). Consider the partial derivatives of $Q_n(\boldsymbol{\beta} + \mathbf{v})$. For $j = k_n + 1, \ldots, p_n$, if $|v_j| \leq \lambda_n$,

$$\begin{aligned}\frac{\partial Q_n(\boldsymbol{\beta}+\mathbf{v})}{\partial v_j} &= 2\sum_{i=1}^n X_{ij}(\varepsilon_i - \mathbf{X}_i'\mathbf{v}) + n\lambda_n \mathrm{sgn}(v_j) \\ &= 2\sum_{i=1}^n X_{ij}\varepsilon_i - 2\sum_{i=1}^n X_{ij}\mathbf{X}_{i1}'\mathbf{v}_1 - 2\sum_{i=1}^n X_{ij}\mathbf{X}_{i2}'\mathbf{v}_2 + n\lambda_n \mathrm{sgn}(v_j) \\ &\triangleq II_{n1,j} + II_{n2,j} + II_{n3,j} + II_{n4,j}.\end{aligned}$$

Examine the first three terms one by one.

$$E[\max_{k_n+1 \leq j \leq p_n} |II_{n1,j}|] \leq E^{1/2}\left[\sum_{j=k_n+1}^{p_n} II_{n1,j}^2\right] = 2\sqrt{nm_n}\sigma,$$

$$\begin{aligned}\max_{k_n+1 \leq j \leq p_n} |II_{n2,j}| &= 2 \max_{k_n+1 \leq j \leq p_n} \left|\sum_{i=1}^n X_{ij}\mathbf{X}_{i1}'\mathbf{v}_1\right| \\ &\leq 2\|\mathbf{v}_1\| \max_{k_n+1 \leq j \leq p_n} \sqrt{(\mathbb{X}_{\cdot j})'\mathbb{X}_1 \mathbb{X}_1' \mathbb{X}_{\cdot j}} \\ &\leq 2\|\mathbf{v}_1\| \max_{k_n+1 \leq j \leq p_n} \|\mathbb{X}_{\cdot j}\| \rho_{\max}^{1/2}(\mathbb{X}_1 \mathbb{X}_1') \\ &= 2\|\mathbf{v}_1\| \max_{k_n+1 \leq j \leq p_n} \|\mathbb{X}_{\cdot j}\| \rho_{\max}^{1/2}(\mathbb{X}_1' \mathbb{X}_1) \\ &= 2n\sqrt{\pi_{n,k_n}} \|\mathbf{v}_1\|,\end{aligned}$$

$$\begin{aligned}\max_{k_n+1 \leq j \leq p_n} |II_{n3,j}| &= 2 \max_{k_n+1 \leq j \leq p_n} |\sum_{i=1}^n X_{ij}\mathbf{X}_{i2}'\mathbf{v}_2| \\ &\leq 2\|\mathbf{v}_1\| \|\mathbb{X}_{\cdot j}\| \rho_{\max}^{1/2}(\mathbb{X}_2' \mathbb{X}_2) \\ &= 2n\sqrt{\omega_{n,m_n}} \|\mathbf{v}_2\|.\end{aligned}$$

Following the above argument we have

$$\begin{aligned}&P\left(\bigcup_{k_n+1 \leq j \leq p_n} \{|II_{n1,j}| > |II_{n4,j}| - |II_{n2,j}| - |II_{n3,j}|\}\right) \\ &\leq \frac{2\sqrt{nm_n}\sigma_2}{n\lambda_n - 2n\left(\sqrt{\pi_{n,k_n}}\|\mathbf{v}_1\| + \sqrt{\omega_{n,m_n}}\|\mathbf{v}_2\|\right)}.\end{aligned}$$

When (A3) holds, $\sqrt{n}\lambda_n/\sqrt{m_n} \to \infty$. Under (A1)–(A2), $\|\mathbf{v}\| = O_P(\sqrt{p_n/n}/\rho_{n,1})$. Therefore

$$P\left(\bigcup_{k_n+1 \leq j \leq p_n} \{|II_{n1,j}| > |II_{n4,j}| - |II_{n2,j}| - |II_{n3,j}|\}\right) \to 0 \text{ as } n \to \infty.$$

This indicates that with probability tending to 1, for all $j = k_n+1, \ldots, p_n$, the sign of $\frac{\partial Q_n(\boldsymbol{\beta}+\mathbf{v})}{\partial v_j}$ is the same as $v_j$, provided that $|v_j| < \lambda_n$, which further implies that

$$\lim_{n \to \infty} P(\widehat{\boldsymbol{\beta}}_{2n} = \mathbf{0}_{m_n}) = 1. \quad \square$$

*Proof of Theorem 4.* Follow the argument in the proof of Theorem 3. Note that in the random design setting, under (B1.a),

$$\max_{k_n+1 \leq j \leq p_n} |II_{n2,j}| = 2 \max_{k_n+1 \leq j \leq p_n} \left| \sum_{i=1}^n X_{ij} \mathbf{X}'_{i1} \mathbf{v}_1 \right|$$

$$\leq 2\|\mathbf{v}_1\| \max_{k_n+1 \leq j \leq p_n} \sqrt{(\mathbb{X}_{\cdot j})' \mathbb{X}_1 \mathbb{X}'_1 \mathbb{X}_{\cdot j}}$$

$$\leq 2\|\mathbf{v}_1\| \max_{k_n+1 \leq j \leq p_n} \|\mathbb{X}_{\cdot j}\| \rho_{\max}^{1/2}(\mathbb{X}_1 \mathbb{X}'_1)$$

$$\leq 2\|\mathbf{v}_1\| \sqrt{n} M \sqrt{\rho_{\max}(\mathbb{X}'_1 \mathbb{X}_1)}$$

$$\leq 2M\sqrt{n}\|\mathbf{v}_1\| \sqrt{n \left[\rho_{\max}(E[\mathbf{X}_1 \mathbf{X}'_1]) + \|\mathbf{A}_{11}\|\right]}$$

$$\leq 2n\|\mathbf{v}_1\| \sqrt{\pi_{k_n} + E^{1/2}\|\mathbf{A}_{11}\|^2 O_P(1)}$$

$$= 2n\|\mathbf{v}_1\| \sqrt{\pi_{k_n} + O_P(\rho_1) \frac{M_4^{1/2} k_n}{\rho_1 \sqrt{n}}}$$

$$\leq 4n\|\mathbf{v}_1\| \sqrt{\pi_{k_n}} O_P(1)$$

for sufficiently large $n$. Similarly

$$\max_{k_n+1 \leq j \leq p_n} |II_{n3,j}| \leq 4n\|\mathbf{v}_2\| \sqrt{\omega_{m_n}} O_P(1).$$

The rest of the argument is identical to that in the fixed design case and thus omitted here. □

*Proof of Theorem 5.* During the course of proving Lemma 1, we have under (A0)–(A1), $\|\widehat{\boldsymbol{\beta}}_n - \boldsymbol{\beta}\| = O_P(\lambda_n \sqrt{k_n/\rho_{n,1}} + \sqrt{p_n/(n\rho_{n,1})})$. Under (A2), this implies that

$$\|\widehat{\boldsymbol{\beta}}_{1n} - \boldsymbol{\beta}_1\| = o_P(\min_{1 \leq j \leq k_n} |\boldsymbol{\beta}_j|).$$

Also from (A2), $\lambda_n = o(\min_{1 \leq j \leq k_n} |\boldsymbol{\beta}_j|)$. Therefore, with probability tending to 1, all the $\widehat{\beta}_j$ ($1 \leq j \leq k_n$) are bounded away from $[-a\lambda_n, a\lambda_n]$ and so the partial derivatives exist. At the same time, $\widehat{\boldsymbol{\beta}}_{2n} = \mathbf{0}_{m_n}$ with probability tending to 1. Thus with probability tending to 1, the stationarity condition holds for the first $k_n$ components. That is, $\widehat{\boldsymbol{\beta}}_{1n}$ necessarily satisfies the equation

$$\sum_{i=1}^n (Y_i - \mathbf{X}'_{i1}\widehat{\boldsymbol{\beta}}_{1n}) \mathbf{X}_{i1} = 0, \quad \text{i.e.} \quad \sum_{i=1}^n \varepsilon_i \mathbf{X}_{i1} = \sum_{i=1}^n \mathbf{X}_{i1} \mathbf{X}'_{i1}(\widehat{\boldsymbol{\beta}}_{1n} - \boldsymbol{\beta}_1).$$

So the random vector being considered

$$\mathbf{Z}_n \triangleq \sqrt{n} \Sigma_n^{-1/2} \mathbf{A}_n (\widehat{\boldsymbol{\beta}}_{1n} - \boldsymbol{\beta}_1)$$

$$= \sqrt{n} \sum_{i=1}^n \Sigma_n^{-1/2} \mathbf{A}_n (\mathbb{X}'_1 \mathbb{X}_1)^{-1} \mathbf{X}_{i1} \varepsilon_i$$

$$\triangleq n^{-1/2} \sum_{i=1}^n \mathbf{R}_i^{(n)} \varepsilon_i,$$

where $\mathbf{R}_i^{(n)} = \Sigma_n^{-1/2}\mathbf{A}_n(n^{-1}\mathbb{X}_1'\mathbb{X}_1)^{-1}\mathbf{X}_{i1}$. The equality holds with probability tending to 1. $\max_{1\le i\le n} \|\mathbf{R}_i^{(n)}\|/\sqrt{n} \to 0$ is implied by (A4), as can be seen from

$$\frac{\|\mathbf{R}_i^{(n)}\|}{\sqrt{n}} = \frac{\|\Sigma_n^{-1/2}\mathbf{A}_n\left(n^{-1}\mathbb{X}_1'\mathbb{X}_1\right)^{-1}\mathbf{X}_{i1}\|}{\sqrt{n}}$$
$$\le n^{-1/2}\left\|\left(n^{-1}\mathbb{X}_1'\mathbb{X}_1\right)^{-1/2}\mathbf{X}_{i1}\right\|$$
$$\cdot \rho_{\max}^{1/2}\left(\left(n^{-1}\mathbb{X}_1'\mathbb{X}_1\right)^{-1/2}\mathbf{A}_n'\Sigma_n^{-1}\mathbf{A}_n\left(n^{-1}\mathbb{X}_1'\mathbb{X}_1\right)^{-1/2}\right)$$
$$= n^{-1/2}\left\|\left(n^{-1}\mathbb{X}_1'\mathbb{X}_1\right)^{-1/2}\mathbf{X}_{i1}\right\|\rho_{\max}^{1/2}\left(\sigma^{-2}\Sigma_n^{-1/2}\Sigma_n\Sigma_n^{-1/2}\right)$$
$$= \sqrt{\sigma^{-2}\mathbf{X}_{i1}'\left(\sum_{i=1}^n \mathbf{X}_{i1}\mathbf{X}_{i1}'\right)^{-1}\mathbf{X}_{i1}}.$$

Therefore for any $\xi > 0$,

$$\frac{1}{n}\sum_{i=1}^n E\left[\|\mathbf{R}_i^{(n)}\varepsilon_i\|^2 1\{\|\mathbf{R}_i^{(n)}\varepsilon_i\| > \sqrt{n}\xi\}\right]$$
$$= \frac{1}{n}\sum_{i=1}^n \mathbf{R}_i^{(n)\prime}\mathbf{R}_i^{(n)} E\left[\varepsilon_i^2 1\{\|\mathbf{R}_i^{(n)}\varepsilon_i\| > \sqrt{n}\xi\}\right]$$
$$\le \frac{1}{n}\sum_{i=1}^n \mathbf{R}_i^{(n)\prime}\mathbf{R}_i^{(n)} E\left[\varepsilon_i^2 1\{|\varepsilon_i| > \sqrt{n}\xi/\max_{1\le i\le n}\|\mathbf{R}_i^{(n)}\|\}\right]$$
$$= \frac{1}{n}\sum_{i=1}^n \mathbf{R}_i^{(n)\prime}\mathbf{R}_i^{(n)} o(1)$$
$$= \frac{1}{n}\sum_{i=1}^n \mathbf{X}_{i1}'\left(n^{-1}\mathbb{X}_1'\mathbb{X}_1\right)^{-1}\mathbf{A}_n'\Sigma_n^{-1}\mathbf{A}_n\left(n^{-1}\mathbb{X}_1'\mathbb{X}_1\right)^{-1}\mathbf{X}_{i1} o(1)$$
$$= \sum_{i=1}^n \operatorname{tr}\left\{\left(n^{-1}\mathbb{X}_1'\mathbb{X}_1\right)^{-1}\mathbf{A}_n'\Sigma_n^{-1}\mathbf{A}_n\left(n^{-1}\mathbb{X}_1'\mathbb{X}_1\right)^{-1}\frac{\mathbf{X}_{i1}\mathbf{X}_{i1}'}{n}\right\} o(1)$$
$$= \operatorname{tr}\left\{\left(n^{-1}\mathbb{X}_1'\mathbb{X}_1\right)^{-1}\mathbf{A}_n'\Sigma_n^{-1}\mathbf{A}_n\right\} o(1)$$
$$= \operatorname{tr}\left\{\Sigma_n^{-1}\mathbf{A}_n\left(\frac{1}{n}\sum_{i=1}^n \mathbf{X}_{i1}\mathbf{X}_{i1}'\right)^{-1}\mathbf{A}_n'\right\} o(1)$$
$$= o(1)d.$$

So
$$\mathbf{Z}_n \xrightarrow{D} N(\mathbf{0}_d, I_d).$$

follows from the Lindeberg-Feller central limit theorem and $\operatorname{Var}(\mathbf{Z}_n) = I_d$. □

*Proof of Theorem 6.* The vector being considered

$$\frac{1}{\sqrt{n}}\Sigma_n^{-1/2}\mathbf{A}_n E^{-1/2}[\mathbf{X}_{i1}\mathbf{X}_{i1}']\sum_{i=1}^n \mathbf{X}_{i1}\mathbf{X}_{i1}'(\widehat{\boldsymbol{\beta}}_{1n} - \boldsymbol{\beta}_1)$$
$$= \frac{1}{\sqrt{n}}\Sigma_n^{-1/2}\mathbf{A}_n E^{-1/2}[\mathbf{X}_{i1}\mathbf{X}_{i1}']\sum_{i=1}^n \varepsilon_i \mathbf{X}_{i1}$$

with probability tending to 1. Let $\mathbf{Z}_{ni} = \frac{1}{\sqrt{n}}\Sigma_n^{-1/2}\mathbf{A}_n E^{-1/2}[\mathbf{X}_{i1}\mathbf{X}'_{i1}]\mathbf{X}_{i1}\varepsilon_i$, $i = 1,\ldots,n$. $\{\mathbf{Z}_{ni}, n = 1,2,\ldots, i = 1,\ldots,n\}$ form a triangular array and within each row, they are i.i.d random vectors. First,

$$\operatorname{Var}\left(\sum_{i=1}^n \mathbf{Z}_{ni}\right) = \operatorname{Var}\left(\Sigma_n^{-1/2}\mathbf{A}_n E^{-1/2}[\mathbf{X}_{i1}\mathbf{X}'_{i1}]\mathbf{X}_{11}\varepsilon\right) = I_d.$$

Second, under (B1.a),

$$\sum_{i=1}^n E\left[\|\mathbf{Z}_{ni}\|^2 1_{\{\|\mathbf{Z}_{ni}\|>\xi\}}\right] = nE\left[\|\mathbf{Z}_{n1}\|^2 1_{\{\|\mathbf{Z}_{n1}\|>\xi\}}\right]$$
$$\leq nE^{1/2}[\|\mathbf{Z}_{n1}\|^4]P^{1/2}(\|\mathbf{Z}_{n1}\|>\xi) = o(1),$$

since

$$\begin{aligned}
E^{1/2}[\|\mathbf{Z}_{n1}\|^4] &= E^{1/2}[(\mathbf{Z}'_{n1}\mathbf{Z}_{n1})^2] \\
&= \frac{1}{n}E^{1/2}\left[\varepsilon^4\left(\mathbf{X}'_{11}E^{-\frac{1}{2}}[\mathbf{X}_{11}\mathbf{X}'_{11}]\mathbf{A}'_n\Sigma_n^{-1}\mathbf{A}_n E^{-\frac{1}{2}}[\mathbf{X}_{11}\mathbf{X}'_{11}]\mathbf{X}_{11}\right)^2\right] \\
&\leq \frac{1}{n}\sigma_4^{1/2}\rho_{\max}(\mathbf{A}'_n\Sigma_n^{-1}\mathbf{A}_n)\,\rho_{\max}(E^{-1}[\mathbf{X}_{11}\mathbf{X}'_{11}])E^{1/2}\left[(\mathbf{X}'_{11}\mathbf{X}_{11})^2\right] \\
&\leq \frac{1}{n}\sigma_4^{1/2}\rho_{\max}(\Sigma_n^{-1}\mathbf{A}_n\mathbf{A}'_n)\,\rho_1^{-1}E^{1/2}\left[\left(\mathbf{X}_{11}^{(n)\prime}\mathbf{X}_{11}\right)^2\right] \\
&= \frac{1}{n}\sigma_4^{1/2}\sigma^{-2}\rho_1^{-1}E^{1/2}\left\{\left(\sum_{j=1}^{k_n} X_{1j}^2\right)^2\right\} \\
&= O\left(\frac{k_n}{n\rho_1}\right),
\end{aligned}$$

and

$$P^{1/2}(\|\mathbf{Z}_{n1}\|>\xi) \leq \frac{E^{1/2}(\mathbf{Z}'_{n1}\mathbf{Z}_{n1})}{\xi} = \frac{\sqrt{d}}{\sqrt{n}\xi},$$

by the Lindeberg–Feller central limit theorem we have

$$n^{-1/2}\Sigma_n^{-1/2}\mathbf{A}_n E^{-1/2}[\mathbf{X}_{i1}\mathbf{X}'_{i1}]\mathbb{X}'_1\mathbb{X}_1(\widehat{\boldsymbol{\beta}}_{1n} - \boldsymbol{\beta}_1) \xrightarrow{D} N(\mathbf{0}_d, I_d). \qquad \Box$$

**Acknowledgments.** JH is honored and delighted to have the opportunity to contribute to this monograph in celebration of Professor Piet Groeneboom's 65th birthday and his contributions to mathematics and statistics. The authors thank the editors and an anonymous referee for their constructive comments which led to significant improvement of this article.